\documentstyle[12pt,amssymb,amstex]{amsart}
\numberwithin{equation}{section}

\def\ds{\displaystyle}

\def\a{\alpha}
\def\b{\beta}
\def\g{\gamma}

\def\l{\lambda}

\def\s{\sigma} 
\def\t{\tau}

\newtheorem{thm}{Theorem}[section]
\newtheorem{lem}[thm]{Lemma}

\begin{document}

\centerline{\bf Uzbek Mathematical Journal, Vol 4. , 1997 }

\begin{center}
{\LARGE\textbf {About summability of Fourier-Laplace series}}
\end{center}
\begin{center}
Anvarjon Akhmedov  \\
Tashkent State University\\

\end{center}

{ \small {\bf Abstract}}

 {\small { In this paper we study the almost
everywhere convergence of the expansions related to the
self-adjoint extension of the Laplace operator. The sufficient
conditions for summability is obtained.
 For the orders of Riesz means, which greater than
critical index $\frac{N-1}{2}$ we established the estimation for
maximal operator of the Riesz means. Note that when order $\alpha$
of Riesz means is less than critical index then for establish of
the almost everywhere convergence requests to use other methods
form proving negative results. We have constructed different
method of summability of Laplace series, which based on spectral
expansions property of self-adjoint Laplace-Beltrami operator on
the unit sphere.}}



  \section{Introduction}
Let $S^N$ be unit sphere in $R^{N+1}.$ Let us denote by
$\lambda_0, \lambda_1,...$ the distinct eigenvalues of the Laplace
operator $-\Delta_S$, arranged in increasing order. Let $H_k$
denote the eigenspace corresponding to $\lambda_k.$ We call
elements of $H_k$ spherical harmonics of degree $k$. It is well
known (see \cite{ST5}) that $dimH_k=a_k:$
\begin{eqnarray}\label{eq12*}
a_k=\left\{
\begin{array}{ll}

\ds \qquad 1,\ \qquad\qquad\qquad\qquad\qquad if \qquad k=0,\\[2mm]

\ds \qquad N , \qquad \qquad\qquad\qquad\qquad if \qquad k=1, \\[2mm]

\ds \frac{(N+k)!}{N!k!}-\frac{(N+k-2)!}{N!(k-2)!},  \qquad if
\qquad k\geq 2 \\

\end{array}
\right.
\end{eqnarray}

Let $\hat{A}$ is a self-adjoint extension of the Laplace operator
$\Delta_{S}$  in $L_2(S^N)$ and if $E_\lambda$ is the
corresponding spectral resolution, then for all functions $f\in
L_2(S^N)$ we have

$$
\hat{A}f=\int_0^{\infty}\lambda dE_\lambda f.
$$

The operator$\hat{A}$ has in$L_2(S^N)$ a complete orthonormal
 system of eigenfunctions
  $$ \{Y_1^{(k)}(x),Y_2^{(k)}(x),...,Y_{a_k}^{(k)}(x)\}\subset H_k, k=0,1,2,...,
$$
corresponding to the eigenvalues $\{\l_k=k(k+N-1)\}, k=0,1,2,...$.

 It is easy to check
that the operators $E_\lambda$ have the form
\begin{equation}\label{PS}
E_\lambda f(x)=\sum_{\lambda_n<\lambda}Y_n(f,x) ,
\end{equation}
 where
\begin{equation}\label{FC}
 Y_k(f,x)=\sum\limits_{j=1}^{a_k}Y_j^{(k)}(x)\int\limits_{S^N}f(y)Y_j^{(k)}(y)d\sigma(y)
 \end{equation}.

  The Riesz means of the partial sums (\ref{PS}) is defined by

\begin{equation}\label{RM}
E_n^\alpha
f(x)=\sum\limits_{k=0}^n\left(1-\frac{\lambda_k}{\lambda_n}\right)^\alpha
Y_k(f,x)
\end{equation}

The most convenient object for a detailed investigation are the
expansions of the form  (\ref{RM}). The  integral (\ref{RM}) may
be transformed writing instead of $\hat{A}$ the integral to the
right in (\ref{FC}) and then changing the order of integration.
This yields the formula
 \begin{equation} \label{RM1}
 E_\lambda^s
 f(x)=\int\limits_{S^N}\Theta^s(x,y,\lambda)f(y)d\sigma(y)
 \end{equation}
with

\begin{equation}\label{05}
\Theta^\alpha(x,y,n)=\sum\limits_{k=0}^n\left(1-\frac{\lambda_k}{\lambda_n}\right)^\alpha
Z_k(x,y).
\end{equation}

For $\alpha=0$ this kernel is called the spectral function of the
Laplace operator for the entire space $S^N$.

The behavior of the spectral expansion corresponding to the the
Laplace-Beltrami operator is closely connected with the
asymptotical behavior of the kernel $\Theta^\alpha(x,y,n).$

In the study of questions of a.e. convergence it is convenient to
introduce the maximal operator
$$
E_*^\alpha f(x)=\sup\limits_{\lambda\geq 0}|E_\lambda^\alpha
f(x)|.
$$

 The basic results of this paper is

\begin{thm}\label{Main}
Let $f\in L_p(S^N),1\leq p\leq 2$ then Riesz means  of order
$s>(N-1)\left(\frac{1}{p}-\frac{1}{2}\right)$ of the
Fourier-Laplace series of the function $f$ converges almost
everywhere on $S^N$ to the $f.$
\end{thm}

\section{Preliminaries}

In this section we are going to prove estimates for maximal
operator. Let us recall some more general definition of harmonic
analysis.

For any two points $x$ and $y$ from $S^N$ we shall denote by
$\g(x,y)$ spherical distance between these two points. Actually,
$\g(x,y)$ is a measure of angle between $x$ and $y.$ It is
obvious, that $\g(x,y)\leq \pi.$

 Spherical ball $B(x,r)$ of radius
$r$ and with the center at point $x$ defined by $B(x,r)=\{y\in
S^N: \g(x,y)<r\}.$ For integrable function $f(x)$ the maximal
function of Hardy-Littlwood
\begin{equation}\label{max}
f^*(x)=\sup_{r>0}\frac{1}{|B(x,r)|}\int\limits_{S^N}|f(y)|d\s(y)
\end{equation}
is finite almost everywhere on sphere. The maximal function $f^*$
plays a major role in analysis and has been much studied
(see.\cite{ST5}). In particular, for any $p>1$ and if $f\in L_p$,
then there exists constant $c_p$, such that
$$
\|f^*\|_{L_p}\leq\ \frac{c_p(N)}{p-1}\|f\|_{L_p},
$$
where $c_p$ has no singularities at point $p=1.$

\begin{thm}\label{L1}
Let $\alpha>\frac{N-1}{2}$ then for all $f\in L_1(S^N)$ we have
\begin{equation}\label{L1}
E_*^\a f(x)\leq\frac{c_\a(N)}{\a-\frac{N-1}{2}}
\left(f^*(x)+f^*(\bar{x})\right)
\end{equation}
where $\g(x,\bar{x})=\pi$.
\end{thm}

The prove can be found in \cite{A1}.

 From the boundness of the maximal
function we obtain
\begin{thm}\label{L12}
Let $\alpha>\frac{N-1}{2}$ then for all $f\in L_p(S^N),p>1$ we
have
\begin{equation}\label{L12}
\|E_*^\a f(x)\|_p\leq\frac{c_\a(N)}{\a-\frac{N-1}{2}}\|f\|_p
\end{equation}
\end{thm}

The statement of the Theorem \ref{L12} we will be using when $p$
approaches 1. For using the Stein's interpolation theorem we have
to set the analogous estimation on the case $p=2.$

Let us for $f\in L_2(S^N)$ and $\a>-1/2$ denote an operator $M^\a$
, which plays main role in the estimate of the Riesz means
\begin{equation}
M^\a f(x)=\sup_{n\geq 1}
\left(\frac{1}{n}\sum\limits_{k=0}^n|E_k^\a f(x)\right).
\end{equation}
If we consider the Riesz means of order $\a+\b$ then we may
connect this means with the operator $M^\b$ as follow:
\begin{lem}\label{L3}
Let $\a>-1/2,\ \b>1/2$ then for all $f\in L_2(S^N)$ we have
\begin{equation}\label{EM}
E_*^{\a+\b} f(x)\leq c_{\a,\b}M^\b f(x)
\end{equation}
\end{lem}
Proof. Let $\a>-1/2,\ \b>1/2$ and $f\in L_2(S^N)$, then
\begin{equation}
E_n^{\a+\b}f(x)=\sum\limits_{k=0}^n\left(1-\frac{\l_k}{\l_n}\right)^{\a+\b}Y_k(f,x).
\end{equation}
Using analogous formula of integration by parts to partial sums we
have
\begin{equation}
E_n^{\a+\b}f(x)=\sum\limits_{k=0}^n\left(\left(1-\frac{\l_k}{\l_n}\right)^\b+
\left(1-\frac{\l_k}{\l_n}\right)^\b\right) E_k^\a f(x).
\end{equation}
Applying the Cauchy's inequality we have
\begin{equation}
|E_n^{\a+\b}f(x)|\leq\left(\sum\limits_{k=0}^n|(\left(1-\frac{\l_k}{\l_n}\right)^\b-
\left(1-\frac{\l_k}{\l_n}\right)^\b|^2\right)^{1/2}
\left(\sum\limits_{k=0}^n|E_k^\a f(x)|^2\right)^{1/2}.
\end{equation}

In view that
\begin{equation}
\left(n\sum\limits_{k=0}^n|\left(1-\frac{\l_k}{\l_n}\right)^\b-
\left(1-\frac{\l_{k+1}}{\l_n}\right)^\b|^2\right)^{1/2}
\leq\frac{1}{2}B(2\b-1,\frac{3}{2}),
\end{equation}
where $B(x,y)$ is Beta function:
\begin{equation}
B(x,y)=\int\limits_0^1t^{x-1}(1-t)^{y-1}dt.
\end{equation}

we have
\begin{equation}
E_*^{s+\a}f(x)\leq C_{s,\a}M^\a f(x)
\end{equation}

So it is not hard to see that, if we obtain estimate for $M^\a$,
then the same estimate true for $E_*^{s+\a}.$ For estimate $M^\a$
let enter a new function $G^\a$, which defined as follow:

$$
 G^\a f(x)=\left(\sum\limits_{n=1}^\infty
\frac{1}{n}|E_n^{\a+1}f(x)-E_n^\a f(x)|^2\right)^{1/2}.
$$
for all $f\in L_2(S^N).$
\begin{lem}\label{L3}
Let $\a>-1/2,$  then for all $f\in L_2(S^N)$ we have
\begin{equation}\label{Mm}
\|G^\a (f)\|_{L_2(S^N)}\leq const\ \|f\|_{L_2(S^N)}
\end{equation}
\end{lem}

Proof. Using orthonormality of the functions $\{Y_k(x)\}$ and
Fubini's theorem about the changeable of integration order, we may
estimate the norm of the function $G^\a$ as follow

$$
\|G^\a(f)\|_{L_2(S^N)}\leq\sum\limits_{n=1}^\infty\frac{1}{n}\sum\limits_{k=1}^n
\left(1-\frac{\l_k}{\l_n}\right)^{2\a}\frac{\l_k^2}{\l_n^2}|Y_k(f,x)|^2.
$$

Therefore, in view
$$
\frac{1}{n}\sum\limits_{k=1}^n
\left(1-\frac{\l_k}{\l_n}\right)^{2\a}\frac{\l_k^2}{\l_n^2}\leq
\frac{1}{2}B(2\a+1,5/2),
$$
we have
$$
\|G^\a(f)\|_{L_2(S^N)}\leq\|f\|_{L_2(S^N)}.
$$

Lemma \ref{L3} proved.

\begin{lem}\label{L4}
Let $\a>-1/2,m=1,2,...$  then for all $f\in L_2(S^N)$ we have

\begin{equation}\label{Mm}
 M^\a(f)\leq M^{\a+m}+G^\a(f)+G^{\a+1}(f)+...+G^{\a+m-1}(f)
\end{equation}
\end{lem}

 Proof. We prove inequality (\ref{Mm}) by the Induction method.
 Let first $m=1$, then we have to prove inequality:

\begin{equation}\label{L4.1}
 M^\a(f)\leq M^{\a+1}+G^\a(f)
\end{equation}

For proving the inequality (\ref{L4.1}) the form of $G^\a$
estimate form below as follow:
$$
 [G^\a f(x)]^2=\sum\limits_{n=1}^\infty
\frac{1}{n}|E_n^{\a+1}f(x)-E_n^\a f(x)|^2 \geq \sum\limits_{k=1}^n
\frac{1}{k}|E_k^{\a+1}f(x)-E_k^\a f(x)|^2\geq
$$
$$
\geq\frac{1}{n}\sum\limits_{k=1}^n |E_k^{\a+1}f(x)-E_k^\a
f(x)|^2\geq\frac{1}{n}\sum\limits_{k=1}^n
||E_k^{\a+1}f(x)|-|E_k^\a f(x)||^2=
$$
$$
=\frac{1}{n}\sum\limits_{k=1}^n
|E_k^{\a+1}f(x)|^2+\frac{1}{n}\sum\limits_{k=1}^n |E_k^\a f(x)|^2-
2\frac{1}{n}\sum\limits_{k=1}^n |E_k^{\a+1}f(x)||E_k^\a f(x)|
$$

So using Caushy's inequality
$$
\sum\limits |a_k||b_k|\leq \sqrt{\sum\limits |a_k|^2}
\sqrt{\sum\limits |b_k|^2}
$$

we get
$$
\frac{1}{n}\sum\limits_{k=1}^n
|E_k^{\a+1}f(x)|^2+\frac{1}{n}\sum\limits_{k=1}^n |E_k^\a f(x)|^2-
2\frac{1}{n}\sum\limits_{k=1}^n |E_k^{\a+1}f(x)||E_k^\a f(x)|\geq
$$
$$
\geq \frac{1}{n}\sum\limits_{k=1}^n
|E_k^{\a+1}f(x)|^2+\frac{1}{n}\sum\limits_{k=1}^n |E_k^\a f(x)|^2-
2\frac{1}{n}\sqrt{\sum\limits_{k=1}^n
|E_k^{\a+1}f(x)|^2}\sqrt{\sum\limits_{k=1}^n |E_k^\a f(x)|^2}=
$$
$$
=\left(\sqrt{\frac{1}{n}\sum\limits_{k=1}^n |E_k^\a f(x)|^2}-
\sqrt{\frac{1}{n}\sum\limits_{k=1}^n |E_k^{\a+1}
f(x)|^2}\right)^2.
$$
So finally for $G^\a$ we have
$$
[G^\a (f)]^2\geq \left(\sqrt{\frac{1}{n}\sum\limits_{k=1}^n
|E_k^\a f(x)|^2}- \sqrt{\frac{1}{n}\sum\limits_{k=1}^n |E_k^{\a+1}
f(x)|^2}\right)^2.
$$
Then in view of definition of $M^\a$ we get inequality
(\ref{L4.1}). So inequality (\ref{L4}) proved when $m=1.$ Assume
that (\ref{L4}) true for all $k<m:$
$$
 M^\a(f)\leq M^{\a+k}+G^\a(f)+G^{\a+1}(f)+...+G^{\a+k-1}(f)
$$
Now we have to extend this inequality for $k+1.$ So, as
$$
M^{\a+k}f(x)\leq M^{\a+k+1}+G^{\a+k}
$$
we have
$$
 M^{\a+k}+G^\a(f)+G^{\a+1}(f)+...+G^{\a+k-1}(f)\leq
$$
$$
\leq M^{\a+k+1}+G^\a(f)+G^{\a+1}(f)+...+G^{\a+k-1}(f)+G^{\a+k}(f).
$$
Which proves the assertion of Lemma \ref{L4}.

\begin{lem}\label{L5}
Let $\a>-1/2,$ then for all $f\in L_2(S^N)$ we have
\begin{equation}\label{L5}
\|M^\a(f)\|_{L_2(S^N)} \leq c_\a \|f\|_{L_2(S^N)}.
\end{equation}
\end{lem}
Proof. If we choose the integer $m$ in (\ref{L4}) such that,
$m>\frac{N}{2}$, then we have $\a+m>(N-1)/2$. It is not hard to
show that for all $\a+m>-1/2$
$$
M^{\a+m}f(x)\leq E_*^{\a+m}f(x).
$$
Due to Theorem \ref{L12} we get
$$
\|M^{\a+m}(f)\|_2\leq \|f\|_2.
$$
Then in view Lemma \ref{L3} we have

$$
 \|M^\a(f)\|_{L_2(S^N)}\leq
\|M^{\a+m}(f)\|_{L_2(S^N)}+\|G^\a(f)\|_{L_2(S^N)}+...+\|G^{\a+m-1}(f)\|_{L_2(S^N)}
$$

 So we have
$$
\|M^{\a+m}(f)\|_{L_2(S^N)}\leq \|f\|_{L_2(S^N)}.
$$
Lemma \ref{L5} proved.

Finally for Riesz means we have more important

\begin{thm}\label{T2}
Let $f\in L_2(S^N),$ then for Riesz means of positive order $\a>0$
we have
\begin{equation}\label{T2}
\|E_*^\a(f)\|_{L_2(S^N)}\leq c_\a \|f\|_{L_2(S^N)}.
\end{equation}
\end{thm}
Consequently, for every $f\in L_2(S^N)$ the Riesz means $E_n^\a f$
of any positive order converge almost everywhere on $S^N.$ For
multiple Fourier series this result is due to Mitchell \cite{MIT},
for spectral expansions of the elliptic operator to Peetre
\cite{PEET}. The question of almost everywhere convergence of
$E_n^\a f, f\in L_2(S^N),N>2$, remains open for $\a=0.$

 For multiple Fourier series this result due to Mitchell
(1955). The question of a.e. convergence of $E_\lambda^\alpha
f(x), f\in L_2(S^N),$ remains open for $\alpha=0.$

\section{Interpolation between $L_p,p>1$ and $L_2$}

In this section we are going to prove the almost everywhere
convergence of spectral expansions by Riesz means at index under
the critical line
$\alpha=(N-1)\left(\frac{1}{p}-\frac{1}{2}\right),1\leq p\leq 2$.
Let us now pass to the case $1<p<2$. It follows from (\ref{L12})
that for such values of $p$ the Riesz means converge above the
critical order. This good result for $p$ near 1, but its gets
cruder when $p$ approaches 2, as for $p=2$ a.e. convergence holds
true for all $\alpha>0.$ There arises the natural desire to
interpolate between (\ref{L12}) and (\ref{T2}).

It is intuitively clear that for intermediate $p, 1<p<2,$ an
analogous estimate must hold with an $\a$ that decreases from
$\frac{N-1}{2}$ to 0. Ordinary interpolation allows one to
interpolate the inequalities (\ref{L12}) and  (\ref{T2}) for fixed
operator and it does not provide the possibility to change the
order $\a$ when passing from $L_1$ to $L_2.$ However, if
dependence of the operators under consideration on the parameter
$\a$ is analytic then one can carry out the interpolation in $\a.$

 Let us state Stein's interpolation theorem in a form
suitable for our purposes.\\

We say that a function $\phi(\tau),\tau\in R,$ has admissible
growth if there exist constants $a<\pi$ and $b>0$ such that

\begin{equation}\label{GR}
|\phi(z)|\leq exp(bexpa|\tau|).
\end{equation}

Let $A_z$ be a family of operators defined for simple functions
(i.e. functions which are finite linear combinations of
characteristic functions of measurable subsets of $S^N$.) We term
the family $A_z$ admissible if for any two simple functions $f$
and $g$ the function
$$
\phi(z)=\int\limits_{S^N}f(x)A_zg(z)dx
$$
is analytic in the strip $0\leq Rez\leq 1$ and has admissible
growth in $Imz,$ uniformly in $Rez$ (this means that we have an
estimate in $Imz$ which is analogous to (\ref{GR}), with constants
$a$ and $b$ independent of $Rez$).

\begin{thm}\label{Stein}
Let $A_z $ be an admissible family of linear operators such that
$$
\|A_{i\tau}\|_{L_{p_0}(S^N)}\leq M_0(\tau)\|f\|_{L_{p_0}(S^N)}, \
1\leq p_0\leq\infty,
$$
$$
\|A_{1+i\tau}\|_{L_{p_1}(S^N)}\leq M_1(\tau)\|f\|_{L_{p_1}(S^N)},
\ 1\leq p_1\leq\infty,
$$
 for all simple functions
$f$ and with $M_j(\tau)$ independent of $\tau$ and admissible
growth. Then there exists for each $t, \ 0\leq t\leq 1,$ a
constant $M_t$ such that for every simple function $f$ holds
$$
\|A_{t}\|_{L_{p_t}(S^N)}\leq M_t(\tau)\|f\|_{L_{p_t}(S^N)}, \
\frac{1}{p_t}=\frac{1-t}{p_0}+\frac{t}{p_1}.
$$
\end{thm}

Some of the most useful objects to which this theorem can be
applied are the Riesz means, which analytically depend on $s.$ In
this case, admissible growth in practise does not cause any great
difficulty, as practically all functions encountered in the
applications have exponential growth. the restriction of the
domain of definition of $A_z$ to simple functions does not
diminish the possibility of interpolation, as the simple functions
constitute a dense subset in $L_p(S^N).$\\

Let us turn attention to the difficulties which arise then. The
interpolation theorem \ref{Stein}, which we would like to apply,
pertains to an analytic family of linear operators, but the
maximal operator $E_*^\a$ is nonlinear. This difficulty is
resolved in the following standard way. Denote by $\Psi$ the class
of positive measurable functions on $S^N$ taking finitely many
different values. If $\mu\in\Psi$ then by the definition of the
maximal operator we have
\begin{equation}\label{M1}
|E_{\mu(x)}^\a f(x)|\leq E_*^\a f(x).
\end{equation}
It is clear that one can pick a sequence $\mu_1\leq\mu_2\leq...$
of elements in $\Psi$ such that
$$
\lim_{k\rightarrow\infty}|E_{\mu(x)}^\a f(x)|= E_*^\a f(x).
$$
This allows us to invert (\ref{M1}) as follows:
\begin{equation}\label{M2}
\sup_{\mu\in\Psi}\|E_{\mu(x)}^\a f(x)\|_{L_p(S^N)}= \|E_{*}^\a
f(x)\|_{L_p(S^N)}.
\end{equation}

Fix now $\mu\in\Psi$ and consider the linear operator $E_{*}^\a$,
which depends analytically on the parameter $\a.$

Let us fix an arbitrary $\varepsilon>0$ and set
$\a(z)=\frac{N-1}{2}z+\varepsilon$. Then the operators
$E_{\mu(x)}^{\a(z)}$ satisfy all conditions of interpolation
theorem of Stein with $p_1>1$ and $p_0=2$

\begin{equation}
 \|E_{\mu(x)}^{\a(1+i\t)}(f)\|_{L_{p_1}(S^N)}\leq A_1
\|f\|_{L_{p_1}(S^N)},
\end{equation}
$$
 A_1=\frac{Ce^{\pi|\t| /2}}{(p_1-1)^2},\ Re\
\a(1+i\t)>\frac{N-1}{2},p_1>1,
$$
and
 \begin{equation}
 \|E_{\mu(x)}^{\a(i\t)}(f)\|_{L_2(S^N)}\leq
A_0 \|f\|_{L_2(S^N)},\ Re\ \a(i\t)>0.
\end{equation}

Consequently for all $t:0<t<1$ we have
\begin{equation}
\|E_{\mu(x)}^{\a(t)}(f)\|_{L_p(S^N)}\leq A_t \|f\|_{L_p(S^N)}.
\end{equation}
  where
$\a(t)=\frac{N-1}{2}t+\varepsilon,\
\frac{1}{p}=\frac{1-t}{2}+\frac{t}{p_1},$ and excluding $t$, we
have, $\a>\frac{N-1}{2}\left(\frac{1}{p}-\frac{1}{2}\right)$  and
 $A_t$:

$$
A_t\leq \frac{const}{(p-1)^2}.
$$

Such that we have proved following

\begin{thm}
Let $f(x)\in L_p(S^N), p>1.$ If the order$\a$ of Riesz means
$E_\l^\a f(x)$ is greater than critical index $\frac{N-1}{2}$,
then for maximal operator $E_*^\a$ we have
\begin{equation}\label{SUB}
\|E_{*}^{\a}(f)\|_{L_p(S^N)}\leq \frac{c_{p,\a}}{(p-1)^2}
\|f\|_{L_p(S^N)},
\end{equation}
where constant $c_{p,\a}$ has no singularities at $p=1$ and
$\a=\frac{N-1}{2}.$
\end{thm}

The estimate (\ref{SUB}) says that the means $E_\lambda^\alpha
f(x)$ of any function $f\in L_p(S^N), 1<p<2,$ converge a.e. for
the values of $\alpha$ indicated. Let us note that the convergence
theorem is valid also for $p=1.$ In this case, the maximal
function $f^*$ is not of strong type (1,1), so the estimate
(\ref{L12}) is not true, but it is of weak type (1,1), which in
view (\ref{L11}) is sufficient for the a.e. convergence of the
means $E_\lambda^\alpha f(x)$ for $\alpha>\frac{N-1}{2}$ for any
integrable function $f.$

\section{Conclusion}

To sum up, we may conclude that a sufficient condition for the
a.e. convergence of the Riesz means $E_\lambda^\alpha f(x)$ of a
function $f\in L_p(S^N)$ is that we have
\begin{equation}\label{COND}
\alpha>(N-1)\left(\frac{1}{p}-\frac{1}{2}\right), 1\leq p <2.
\end{equation}
How sharp is this condition? To elucidate this question we begin
with the case of Riesz means of order $\alpha=\frac{N-1}{2}.$
Condition (\ref{COND}) shows that for the a.e. convergence of
Riesz means one requires an order above the critical index only if
$p=1.$ That this requirement is essential follows from the
following theorem:
\begin{thm}
There exists a function $f\in L_1(S^N),N\geq 2,$ such that almost
everywhere on $S^N$
$$
\lim\limits_{\lambda\rightarrow\infty}|E_\lambda^\alpha
f(x)|=+\infty, \alpha=\frac{N-1}{2}.
$$
\end{thm}
Thus, for $p=1$ condition (\ref{COND}) is sharp. For $p>1$ we have
(see \cite{A2}):

\begin{thm}
If
\begin{equation}\label{COND2}
 0\leq \alpha\leq
N\left(\frac{1}{p}-\frac{1}{2}\right)-\frac{1}{2}, 1\leq p\leq
\frac{2N}{N+1},
\end{equation}
then there exists a function $f\in L_p(S^N)$ such that
$E_\lambda^\alpha f(x)$ is divergent on a set of positive measure.
\end{thm}

Between the conditions (\ref{COND}) and (\ref{COND2}) there is a
gap.

\end{document}